\documentclass[11pt,letterpaper]{amsart}
\usepackage[all]{xy}
\usepackage{amssymb}
\usepackage{amsthm} 
\usepackage{cite}
\usepackage{hyperref}
\usepackage[normalem]{ulem}
\usepackage{amsmath,mathtools}
\usepackage{amscd,enumitem}
\usepackage{graphicx}
\usepackage{caption}
\usepackage[justification=centering]{caption}
\usepackage{verbatim}
\usepackage{eurosym}
\usepackage{float}
\usepackage{color}
\usepackage{dcolumn}
\usepackage[mathscr]{eucal}
\usepackage[all]{xy}
\usepackage{bbm}
\usepackage{ytableau}
\newtheorem*{conj*}{Conjecture}
\newtheorem{theorem}{Theorem}[section]

\theoremstyle{definition}
\newtheorem*{remark}{Remark}

\theoremstyle{plain}

\newtheorem{lemma}[theorem]{Lemma}

\newtheorem{prop}[theorem]{Proposition}

\newtheorem{definition}[theorem]{Definition}

\newcommand{\Z}{\mathbb{Z}}

\newcommand{\Q}{\mathbb{Q}}

\newcommand{\N}{\mathbb{N}}

\newcommand{\C}{\mathbb{C}}

\newcommand{\Zhat}{\widehat{Z}}
\newcommand{\Zhathat}{\widehat{\vphantom{\rule{5pt}{10pt}}\smash{\widehat{Z}}\,}\!}

\numberwithin{equation}{section}

\newtheoremstyle{example}
  {\topsep}   
  {\topsep}   
  {\normalfont}  
  {0pt}       
  {\bfseries} 
  {.}         
  {5pt plus 1pt minus 1pt} 
  {}          
\theoremstyle{example}

\def\({\left(}
\def\){\right)}

\usepackage{centernot}

\begin{document}
\title{Infinite Families of Quantum Modular 3-Manifold Invariants}
\author{Louisa Liles and Eleanor McSpirit}

\address{Department of Mathematics, University of Virginia, Charlottesville VA 22904}
\email{lml2tb@virginia.edu, egm3zq@virginia.edu}

\thanks{LL acknowledges the partial support of NSF RTG Grant DMS-1839968, NSF Grant DMS-2105467, and the UVa Grad Council Research Grant. EM acknowledges the partial support of the UVa Dean's Doctoral Fellowship, as well as NSF Grants DMS-9508976 and DMS-9304580.}

\begin{abstract}
One of the first key examples of a quantum modular form, which unifies the Witten-Reshetikhin-Turaev (WRT) invariants of the Poincar\'e homology sphere, appears in work of Lawrence and Zagier. We show that the series they construct is one instance in an infinite family of quantum modular invariants of negative definite plumbed 3-manifolds whose radial limits toward roots of unity may be thought of as a deformation of the WRT invariants. We use a recently developed theory of Akhmechet, Johnson, and Krushkal (AJK) which extends lattice cohomology and BPS $q$-series of 3-manifolds. As part of this work, we provide the first calculation of the AJK series for an infinite family of $3$-manifolds. Additionally, we introduce a separate but related infinite family of invariants which also exhibit quantum modularity properties. 
\end{abstract}

\maketitle

\section{Introduction and Statement of Results}
In 1999, Lawrence and Zagier established a connection between modular forms and invariants of 3-manifolds arising from quantum topology and physics \cite{lawrencezagier}. The Witten-Reshetihkin-Turaev (WRT) invariants, conceptualized by Witten in terms of a path integral for $SU(2)$ Chern-Simons theory and made mathematically precise by Reshetikhin and Turaev, are a family of 3-manifold invariants indexed by roots of unity \cite{witten, R-T}.

Lawrence and Zagier unified these WRT invariants by defining a holomorphic function $A(q)$ on the unit disk whose limiting values at roots of unity recover the invariants of the Poincar\'e homology sphere. The function $A(q)$ is in some sense an integral of a theta function; its obstruction to modularity on the unit circle is given by a factor introduced by integrating. It is one of the first key examples of a \textit{quantum modular form}, a term coined by Zagier in 2010 with this series in mind \cite{zagier}. Zagier's seminal work has inspired an extensive and ongoing body of research on quantum modular forms; see Chapter 21 of \cite{harmonic} and the references listed therein.

Gukov, Pei, Putrov, and Vafa (GPPV) provided a physical interpretation for $A(q)$ which not only gave rise to the construction of analogous series $\Zhat(q)$ for a large class of $3$-manifolds  equipped with a spin$^{c}$ structure, but also motivated a conjecture that such series exist for all $3$-manifolds \cite{gppv,gukov_2017}. These unified invariants were further extended in \cite{gukovmanolescu, gukovmodularity} and the quantum modularity of these series for certain classes of manifolds were established in \cite{bringmannquantum1, bringmannquantum2}.

Recent work of Akhmechet, Johnson, and Krushkal generalizes $\Zhat(q)$ to a two-variable series invariant $\Zhathat_{Y}(t,q)$ of a $3$-manifold $Y$ \cite{latticecohomology}. This series is defined using an extension of lattice cohomology, a theory developed by N\'emethi in \cite{nemethi}. More details on the topological construction and properties of this invariant appear in Section \ref{sec: lattice cohomology}.

The aim of this paper is to show that this two-variable series gives rise to infinitely many quantum modular forms whose values at roots of unity can be considered deformations of the WRT invariants. These results realize the work of Lawrence and Zagier as a special case. As a first step in this process, we provide a calculation of $\Zhathat_{\Sigma}(t,q)$ where $\Sigma$ is a Brieskorn homology sphere, generating the first known calculation of this invariant for an infinite family of manifolds.

In particular, we find explicit formulae for the coefficients $\varphi(n;t)$ of the $q$-series $\Zhathat$, which are Laurent polynomials in $t$. The result is, for each Brieskorn sphere $\Sigma$,
a $q$-series of the form
\begin{equation}\label{eqn:qseries} \Zhathat_\Sigma(t,q) = q^\Delta \Big( C - \sum_{n \geq 0} \varphi(n;t)q^\frac{n^2}{4p} \Big),\end{equation}
where $\Delta \in \mathbb{Q}$, \ $p \in \mathbb{Z}$, and $C$ is zero unless $\Sigma$ is the Poincar\'e homology sphere, in which case it equals $q^{1/120}(t+t^{-1})$
; see Section \ref{sec: series analysis} for full definitions.

A priori, $\Zhathat_\Sigma(t,q)$ is convergent as a two-variable series for $t \in \mathbb{C}$ and $|q| < 1$. By leveraging the arithmetic properties of the coefficients $\varphi(n;t)$ when $t$ is a root of unity, we are able to show the following:

\begin{theorem} \label{thm: radial limits} Let $\zeta$ be a $j$th root of unity, $\xi$ a $K$th root of unity, and $\Sigma$ a Brieksorn sphere. Define $\Zhathat_\Sigma(\zeta,\xi) := \displaystyle{\lim_{{t \searrow 0}} \Zhathat_\Sigma(\zeta,\xi e^{-t})}$. This limit exists and we have
\[\Zhathat_\Sigma(\zeta,\xi) = \xi^\Delta \bigg( D + \sum_{n=1}^{2pjK} \left(\frac{n}{2pjK}-\frac{1}{2} \right) \varphi(n;\zeta)\xi^{\frac{n^2}{4p}} \bigg),\]
where $D=\xi^{1/120}\text{Re}(\zeta)$ when $\Sigma$ is the Poincar\'e homology sphere and zero otherwise. 
\end{theorem}

In general, these limit calculations give a family of ``$t$-deformed'' WRT invariants whose topological interpretation is an open question. However, using the above results we prove that for $t$ a fixed root of unity, $\Zhathat_{\Sigma}(t,q)$ is, up to normalization, a quantum modular form:

\begin{theorem} \label{thm:modularity} Let $q = e^{2\pi i \tau}$. If $\zeta$ is a $j$th root of unity and $\Sigma$ is a Brieskorn sphere, then
\[\Zhathat_\Sigma(\zeta,q)  = q^\Delta\left(C - A_\zeta(\tau)\right), \]
where $A_\zeta(\tau)$ is a quantum modular form of weight $1/2$ with respect to $\Gamma(4pj^2)$.
\end{theorem}

\begin{remark}
 The definition of a quantum modular form is deferred until Section \ref{sec: modularity discussion}. 
We further have that $A_\zeta(\tau)$ is a ``strong" quantum modular form in the sense of \cite{zagier}.
\end{remark}

The classical theory of theta functions involves forms of weight $1/2$ and $3/2$ that are related through differentiation of the Jacobi theta function.
Because of the existence of the second variable, one can differentiate $\Zhathat_\Sigma(t,q)$, summand by summand, with respect to $t$ and consider the new invariant that arises. This series, under specialization, also enjoys quantum modularity properties; the result is a sum of quantum modular forms of mixed weight:

\begin{theorem}\label{thm:diffmodular} Define $\displaystyle{{\Zhathat'}_{\Sigma}(t, q):= t \frac{\partial}{\partial t} {\Zhathat}_{\Sigma}(t, q)}$. Let $\zeta$ be a $j$th root of unity, and let $C'$ equal $q^{1/120}(t-t^{-1})$ when $\Sigma$ is the Poincar\'e homology sphere and equal zero otherwise. Then \[{\Zhathat'}_{\Sigma}(\zeta, q)= q^\Delta\left(C' - A'_\zeta(\tau)\right), \] where $A'_\zeta(\tau)$ is a sum of quantum modular forms of weight $1/2$ and $3/2$ for $\Gamma(4pj^2)$. 
\end{theorem}
\begin{remark}
The \textit{quantum set} is notably smaller for the quantum modular forms in Theorem \ref{thm:diffmodular}. This is due to the fact that the weight $3/2$ quantum modular form need not correspond to the Eichler integral of a cusp form. For more details, see Section \ref{sec: modularity proof} and Theorem 1.1 of \cite{goswamiosburn}.
\end{remark}

The paper is organized as follows. In Section \ref{sec: lattice cohomology} we recall the necessary background from low-dimensional topology to understand and motivate the study of this two-variable series. In Section \ref{sec: series analysis} we derive an explicit formula for the $q$-series coefficients of $\Zhathat_{\Sigma}(t,q)$ where $\Sigma$ is a Brieskorn sphere. Section \ref{sec: radial limits} contains a proof of Theorem \ref{thm: radial limits}. In Section \ref{sec: modularity discussion} we discuss the theory of modular and quantum modular forms as it pertains to $\Zhathat_{\Sigma}(t,q)$, and in Section \ref{sec: modularity proof} we offer proofs of Theorems \ref{thm:modularity} and \ref{thm:diffmodular}.

\section*{Acknowledgements}
The authors are thankful for the guidance and thoughtful comments provided by Vyacheslav Krushkal, Ken Ono, Peter Johnson, and Ross Akhmechet.

\section{The $\Zhathat$ Invariant}\label{sec: lattice cohomology}

We begin with a motivating overview of the $\Zhathat$ invariant. This recently-developed two-variable series provides a common refinement of two existing invariants, the GPPV invariant $\Zhat$ and lattice cohomology. To any negative-definite plumbed $3$-manifold equipped with a spin$^{c}$ structure, one can associate a combinatorial object encoding the $0^\text{th}$ lattice cohomology, called a graded root. The authors of \cite{latticecohomology} assigned Laurent polynomial weights to the vertices of this root such that this new ``weighted graded root" is still an invariant of the manifold. The series $\Zhathat_{Y}(t,q)$ results from taking the limit, in a precise sense, of these weights. Setting $t=1$ recovers the GPPV invariant $\Zhat(q)$. This paper calculates the two-variable series $\Zhathat$ for an infinite family of $3$-manifolds, however the calculation of weighted graded roots remains an open problem. Below, we cover details of this construction necessary for our work.

\subsection{Negative Definite Plumbed $3$-Manifolds}
Let $\Gamma$ be a finite graph with integer weights on its vertices.
As in \cite{latticecohomology}, we restrict to the case in which $\Gamma$ is a tree. Let $m:v(\Gamma)\to \Z$ be the corresponding weight function and $s=|v(\Gamma)|$. Choosing an order on $v(\Gamma)$ enables us to write a weight vector $m \in \Z^{s}$ given by $m_{i}=m(v_{i})$ and a degree vector $\delta \in \Z^{s}$ given by $\delta_{i}=\delta(v_{i})$. With this ordering, we can associate to $\Gamma$ a symmetric $s \times s$ matrix $M$ given by \[M_{i,j}=\begin{cases}m_{i} & i=j \\ 1 & i \neq j\text{ and $v_{i}$ and $v_{j}$ are connected by an edge} \\ 0 & \text{otherwise.}\end{cases}\] We say $\Gamma$ is negative definite whenever $M$ is negative definite.

To obtain a $3$-manifold from $\Gamma$, create a framed link $\mathcal{L}(\Gamma) \subset S^{3}$ by associating to each vertex $v_{i}$ an unknot with framing $m_{i}$ and Hopf linking unknots together whenever their corresponding vertices share an edge. The resulting linking matrix of $\mathcal{L}(\Gamma)$ is the plumbing matrix $M$. $Y(\Gamma)$ is defined to be the $3$-manifold obtained by Dehn surgery on $\mathcal{L}(\Gamma)$. Equivalently, $Y=\partial X$ where $X$ is obtained by adding $2$-handles to $\mathbb{D}^{4}$ along $\mathcal{L}(\Gamma)$. From this perspective, $M$ represents the intersection form of $X$. 

\begin{figure}[h]\label{tree2715}
\includegraphics[scale=0.19]{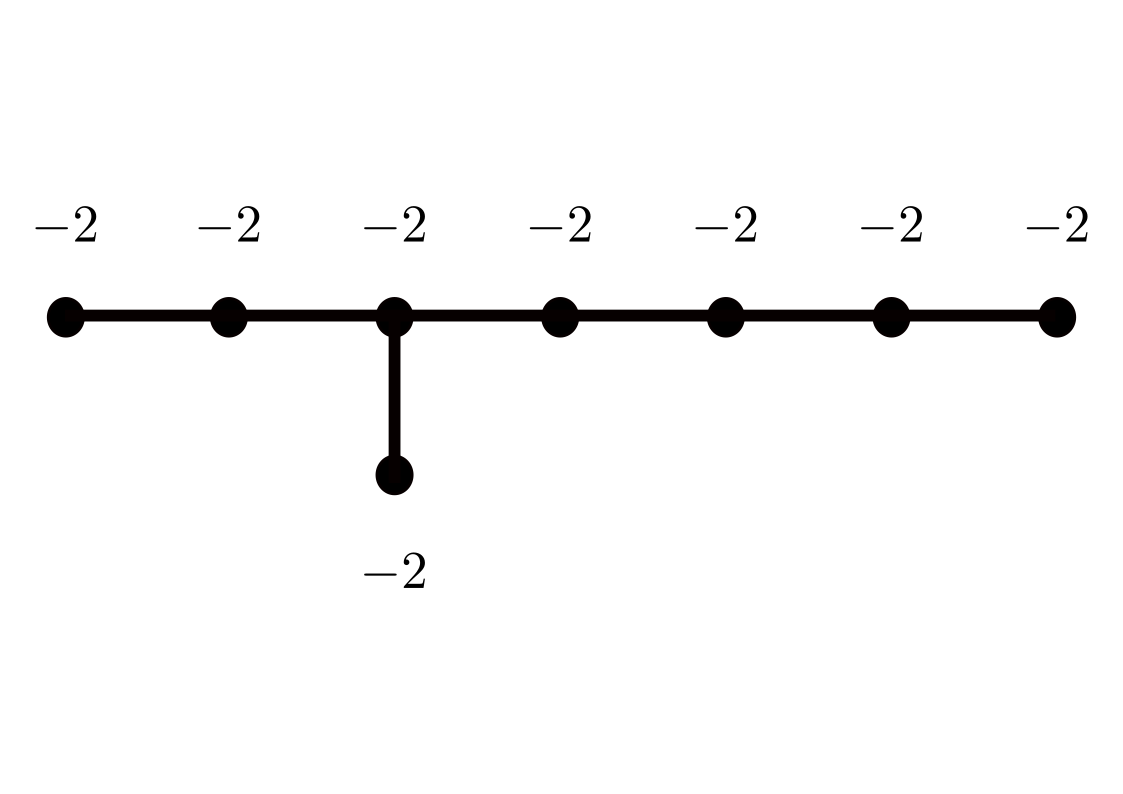} \qquad 
\includegraphics[scale=0.19]{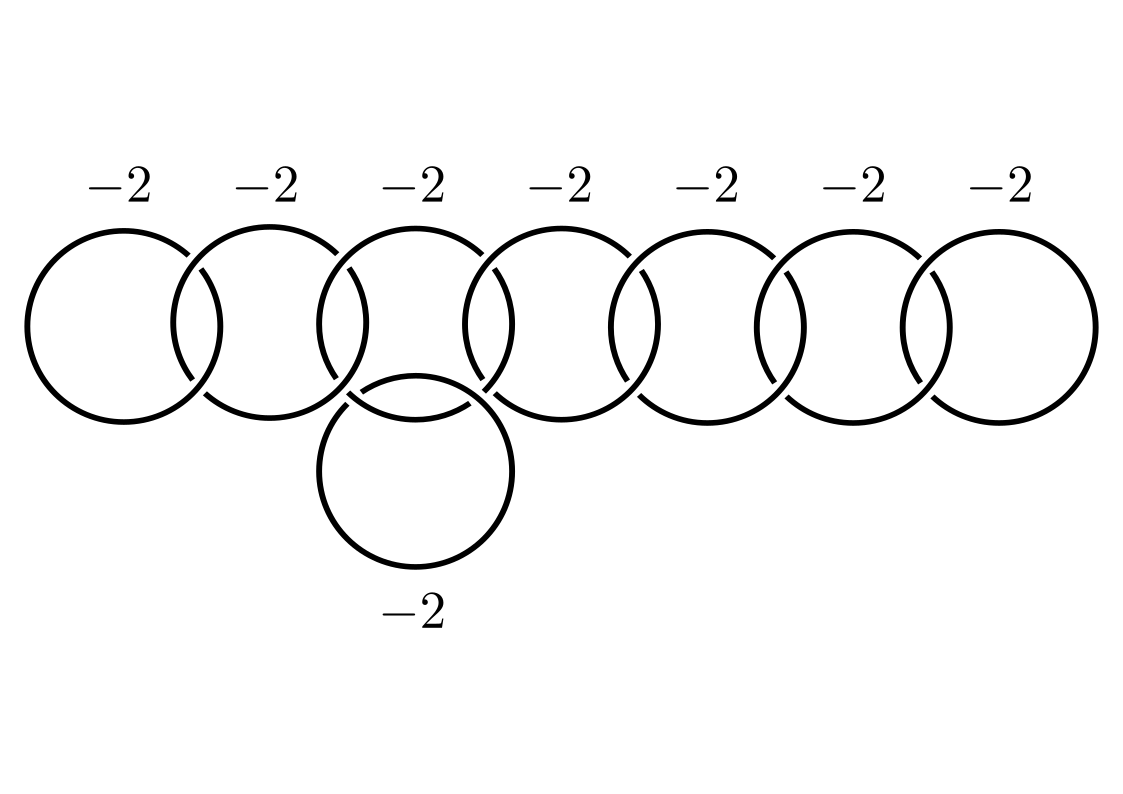} \vspace{-.2in}
\caption{A plumbing tree and its associated link for the Poincar\'e homology sphere $\Sigma(2,3,5)$}
\label{fig:tree2715}
\end{figure}

In general, we say $Y$ is a negative-definite plumbed $3$-manifold if it is homeomorphic to $Y(\Gamma)$ for some negative-definite plumbing graph $\Gamma$. Two distinct plumbing trees may result in homeomorphic manifolds; in fact this is the case if and only if the trees can be related by a finite sequence of Neumann moves of type (a) and (b) \cite{neumann}. Therefore, an invariant of a negative-definite plumbed manifold $Y$ must be invariant under these two moves on its plumbing graph.

As with the GPPV invariant, $\Zhathat_{Y}(t,q)$ takes as inputs a negative-definite plumbed $3$-manifold $Y$ and a chosen spin$^{c}$ structure. The set of spin$^{c}$ structures can be given in terms of the plumbing data; it is known that \begin{equation}\label{eqn:spinc}\text{spin}^{c}(Y)\cong \frac{m+2\Z^{s}}{2M\Z^{s}}\cong \frac{\delta+2\Z^{s}}{2M\Z^{s}}\end{equation} where the second isomorphism is given by $[k] \mapsto [k -(m +\delta)]$. For a more detailed discussion of spin$^{c}$ structures, see Section 2.2 of \cite{latticecohomology}. 
\subsubsection{Key example: Brieskorn homology spheres}  Let $(b_{1},b_{2},b_{3})$ be pairwise relatively prime positive integers $b_{1} < b_{2} < b_{3}$. The corresponding Brieskorn sphere
$\Sigma(b_{1},b_{2},b_{3})$ is given by \[\Sigma(b_{1},b_{2},b_{3})=\{(z_{1},z_{2},z_{3}) \in \C^{3} \ \mid \ z_{1 }^{b_{1}}+z_{2}^{b_{2}}+z_{3}^{b_{3}}=0\}\cap S^{5} \subset \C^{3},\] the intersection of a singular complex surface with the unit sphere in $\C^{3}$.
Given $(b_{1},b_{2},b_{3})$, Neumann and Reymond provide an algorithm by which one can find a plumbing tree $\Gamma$ for the associated Brieskorn sphere \cite{neumannreymond}. This process guarantees that $\Gamma$ is always a star graph with one $3$-valent vertex and $3$ legs, as is the case in Figure \ref{fig:tree2715}. 

As integral homology spheres, Brieskorn spheres have only one spin$^{c}$ structure, and $\Zhathat$ is independent of choice of spin$^{c}$ representative. Therefore, in calculations involving Brieskorn spheres we will drop the subscript indicating the spin$^{c}$ structure.

The general formula for $\Zhathat_{Y}(t,q)$ involves plumbing data, but in Section \ref{sec: series analysis} we give a formula for Brieskorn spheres which only depends on $(b_{1},b_{2},b_{3})$. To achieve this, we use methods similar to those of Gukov and Manolescu, who provide a formula for the GPPV invariant in terms of $(b_{1},b_{2},b_{3})$; see Proposition $4.8$ of \cite{gukovmanolescu}.

\subsection{The Two-Variable Series}
For a choice $k \in \Z^{s}$ of a spin$^{c}$ representative $\displaystyle{[k] \in \frac{m+2\Z^{s}}{2M\Z^{s}}}$, and for any  $x \in \Z^{s}$, let \[\chi_{k}(x):=\frac{-k \cdot x + \langle x, x \rangle}{2} \in \Z,\] where $(\cdot)$ denotes the Euclidean dot product and $\langle -,- \rangle$ denotes the bilinear form given by the plumbing matrix $M$. For $r \in \Z$ and $n \in \mathbb{N}$, let \begin{equation}\label{eq: fhat} \widehat{F}_{n}(r):=\begin{cases} \frac{1}{2}\text{sgn}(r)^{n}\binom{\frac{n+|r|}{2}-2}{n-3} & |r| \geq n-2,r\equiv{n}\bmod{2} \\ 0 & \text{otherwise}, \end{cases}\end{equation} and define \[\widehat{F}_{\Gamma,k}(x):=\prod_{v_{i} \in v(\Gamma)}\widehat{F}_{\delta_{i}}((2Mx+k-Mu)_{i}).\] Note that $\widehat{F}_{n}(r):\Z\to \Q$ describes the coefficient on $z^{-r}$ in the expansion of $(z-z^{-1})^{2-n}$. To state the definition of $\Zhathat$ we will use for our calculations, we define $u:=(1,1,\ldots )$ as well as
\[\Theta_{k}=\frac{k \cdot u - \langle u, u \rangle}{2}\text{,} \qquad \varepsilon_{k}=-\frac{(k-Mu)^{2}+3s+\sum_{v}m_{v}}{4}+2\chi_{k}(x)+\langle x, u \rangle.\]
Then we have the following: 
\begin{theorem}[Theorems $6.3$ and $7.6$ of \cite{latticecohomology}]\label{thm: zhathat def} Let $Y$ be a negative-definite plumbed $3$-manifold with spin$^{c}$ structure $[k]$. The series \begin{equation}\label{zhathat formula} \Zhathat_{Y,[k]}(t,q):=\sum_{x \in \Z^{s}}\widehat{F}_{\Gamma,k}(x)q^{\varepsilon_{k}(x)}t^{\Theta_{k}+\langle x, u \rangle}\end{equation} is an invariant of the pair $(Y, [k])$, and \[\Zhat_{a}(q)=\Zhathat_{Y,[k]}(1,q),\] where $\widehat{Z}_{a}(q)$ is the GPPV invariant for $(Y, [a])$ and $a$ corresponds to $k$ via $(\ref{eqn:spinc})$.
\end{theorem}

\begin{remark} The family of functions $\{\widehat{F}_{n}\}_{n \in \N}$ is defined in \cite{latticecohomology} to be \textit{admissible} in that it satisfies axioms that guarantee that the series analogous to (\ref{zhathat formula}) is an invariant. In this sense $\Zhathat$ belongs to a family of two-variable series invariants developed in \cite{latticecohomology}, parametrized by admissible families. However, if one restricts to Brieskorn spheres, a straightforward calculation shows that any choice of admissible family results in the same two-variable series invariant.
\end{remark}

\section{Series analysis} \label{sec: series analysis}

We now develop an explicit formula for the coefficients of $\Zhathat_{\Sigma}(t,q)$ as a series in $q$ whenever $\Sigma$ is a Brieskorn sphere. The arithmetic properties of these coefficients will allow us to take limits toward roots of unity and establish quantum modularity properties in Sections \ref{sec: radial limits} and \ref{sec: modularity proof}. For a general negative-definite plumbed $3$-manifold $Y$, one can use a program created by Peter Johnson \footnote{Available at \href{https://github.com/peterkj1/plum}{https://github.com/peterkj1/plum}} to calculate the first $N$ coefficients of $\Zhathat_{Y}(t,q)$.

Let $k$ be a spin$^c$ representative for the unique spin$^c$ structure $[k]$ of $\Sigma$, and set $a = k-Mu$. For $x \in \Z^{s}$ we let $\ell:=a+2Mx$. Using the fact established in \cite{latticecohomology} that $\frac{\ell^{T}M\ell}{4}=\frac{a^{2}}{4}-2\chi_{k}(x)-\langle x, u \rangle$, we write

\[\Zhathat_\Sigma(t,q)=q^{-\frac{3s+\sum_{v}m_{v}}{4}}\sum_{x \in \Z^{s}}\prod_{v_{i}\in v(\Gamma)}{\widehat{F}}_{\delta_{i}}(\ell_{i})q^{-\frac{\ell^{T}M^{-1}\ell}{4}}t^{\Theta_{k} + \langle x, u \rangle}.\] Below is the plumbing graph of a Brieskorn sphere with the vertices ordered as needed for this section.
\begin{figure}[h]\label{fig:b graph}
\includegraphics[scale=0.17]{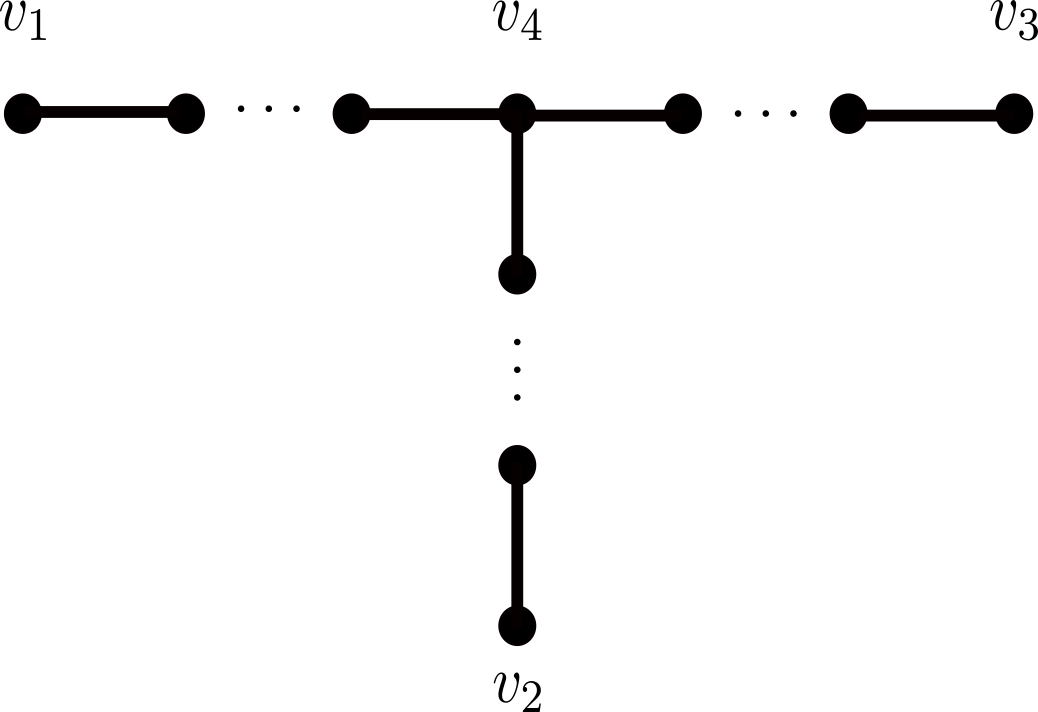}
\caption{the plumbing graph for a Brieskorn sphere $\Sigma$}
\end{figure}
The only $x \in \Z^{s}$ for which $\prod_{v_{i}}{\widehat{F}}_{\delta_{i}}(\ell_i) \neq 0$  are those of the form $\ell=(\varepsilon_{1},\varepsilon_{2},\varepsilon_{3},m, 0, \ldots)$ for $\varepsilon_{i}\in \{\pm 1\}$ and $m$ odd. In this case, we have that $\widehat{F}_{1}(\varepsilon_i)=-\varepsilon_i$ and $\widehat{F}_{3}(m)=\frac{1}{2}\text{sign}(m)$, so\[\prod_{v_{i}\in v(\Gamma)}{\widehat{F}}_{\delta_{i}}(\ell_{i})=-\frac{1}{2}\varepsilon_{1}\varepsilon_{2}\varepsilon_{3}\text{sign}(m).\] Since Brieskorn spheres have unimodular plumbing matrices $M$, every possible combination \\$(\varepsilon_{1},\varepsilon_{2},\varepsilon_{3}, m, 0, \dots)$ is in $a+2M\Z^{s}$. Therefore we can write 
\[\Zhathat_{\Sigma}(t,q)=\frac{-q^{-\frac{3s+\sum_{v}m_{v}}{4}}}{2}\sum_{\varepsilon_{i}\in\{\pm 1\}}\sum_{m\text{ odd}}\varepsilon_{1}\varepsilon_{2}\varepsilon_{3}\text{sign}(m)q^{-\frac{\ell^{T}M^{-1}\ell}{4}}t^{\Theta_{k} + \langle x, u \rangle } \]
One can check that $\langle u, x \rangle = \frac{\varepsilon_{1}+\varepsilon_{2}+\varepsilon_{3}+m-a^{T}u}{2}$. Moreover, since $\Zhathat$ does not depend on a choice of spin$^{c}$ representative, choose $a$ to be $(1,1,1,1,0,\ldots) \in \delta+2M\Z^{s}$. In this case, $\Theta_{k}=2$ and the exponent on $t$ becomes $(\varepsilon_{1}+\varepsilon_{2}+\varepsilon_{3}+m)/2.$

\begin{remark} Following Section 4.6 of \cite{gukovmanolescu} we can can rewrite
\[\frac{-\ell^{t}M^{-1}\ell}{4}=\frac{b_{1}b_{2}b_{3}}{4}\left(m+\sum_{i}\frac{\varepsilon_{i}}{b_{i}}\right)^{2}-\frac{b_{1}b_{2}b_{3}}{4}\sum\frac{1}{b_{i}^{2}}+\frac{\sum_{i}h_{i}}{4},\] 
where $h_{i}$ refers to the cardinality of $H_{1}(\Sigma')$ for $\Sigma'$ the plumbed manifold that results from removing the $i$th vertex of the plumbing graph for $\Sigma$.
Setting 
\[\Delta= \frac{1}{4}\left(\sum_{i}h_{i}-3s-\sum_{v}m_{v} - \frac{b_{2}b_{3}}{b_{1}} - \frac{b_{1}b_{3}}{b_{2}} - \frac{b_{1}b_{2}}{b_{3}}\right)\] and $\varepsilon:= \varepsilon_1+\varepsilon_2+\varepsilon_3$, we now have 
\begin{equation}\label{gm zhathat} \Zhathat_\Sigma(t,q)=\frac{-q^{\Delta}}{2}\sum_{m\text{ odd}}\sum_{\varepsilon_{i}\in\{\pm 1\}}\varepsilon_{1}\varepsilon_{2}\varepsilon_{3}\text{sign}(m)q^{\frac{b_{1}b_{2}b_{3}}{4}(m+\sum_{i}\frac{\varepsilon_{i}}{b_{i}})^{2}}t^{\frac{\varepsilon+m}{2}}.\end{equation}
\end{remark}

Now, set $p:= b_{1}b_{2}b_{3}$ and
\begin{align*}
\alpha_1 &:=b_1b_2b_3 - b_1b_2 -b_1b_3 - b_2b_3; \\
\alpha_2 &:=b_1b_2b_3 + b_1b_2 -b_1b_3 - b_2b_3; \\
\alpha_3 &:=b_1b_2b_3 - b_1b_2 + b_1b_3 - b_2b_3; \\
\alpha_4 &:=b_1b_2b_3 + b_1b_2 + b_1b_3 - b_2b_3.
\end{align*} 

\begin{theorem} \label{thm: q-series} Let $\Sigma(b_{1},b_{2},b_{3})$ be a Brieskorn sphere. Then 
\[\Zhathat_\Sigma(t,q) = q^\Delta \bigg( C - \sum_{n \geq 1} \varphi(n;t)q^\frac{n^2}{4p} \bigg),\]
where $C$ is nonzero and equals $q^{\frac{1}{120}}(t+t^{-1})$ only when $(b_{1},b_{2},b_{3})=(2,3,5)$ and
\[\varphi(n;t) = \begin{cases}
\mp \frac{1}{2}(t^{\frac{\mp n + (\alpha_1+2p)}{2p}} + t^{\frac{\pm n - (\alpha_1+2p)}{2p}}) & n \equiv \pm \alpha_1 \bmod{2p},  \\
\pm \frac{1}{2}(t^\frac{\mp n + \alpha_k}{2p} + t^{\frac{\pm n - \alpha_k}{2p}}) & n \equiv \pm \alpha_k \bmod{2p}, \ k=2,3\\
\mp \frac{1}{2}(t^{\frac{\mp n + (\alpha_4-2p)}{2p}} + t^{\frac{\pm n - (\alpha_4-2p)}{2p}}) & n \equiv \pm \alpha_4 \bmod{2p}, \\ 
0 & otherwise. 
\end{cases}\]
\end{theorem}

\begin{remark}
Note that when $t=1$, this collapses back to the GPPV invariant as calculated by Gukov and Manolescu in \cite{gukovmanolescu}. Fixing $\Sigma=\Sigma(2,3,5)$, the function $\varphi(n;1)$ is equal to $\chi_{+}(n)$ as defined in \cite{lawrencezagier}.
\end{remark}
\begin{proof}
Begin with the calculation given by (\ref{gm zhathat}). Using the fact that \\ $(\varepsilon_{1})(\varepsilon_{2})(\varepsilon_{3})(\text{sign}(m))=(-\varepsilon_{1})(-\varepsilon_{2})(-\varepsilon_{3})(\text{sign}(-m))$ and replacing $m$ odd with $2n+1$, we write \[\Zhathat_\Sigma(t,q)=\frac{-q^{\Delta}}{2}\sum_{\varepsilon_{i}\in \{\pm1\}}\sum_{n \geq 0}\varepsilon_{1}\varepsilon_{2}\varepsilon_{3}q^{p(n^{2}+n+\frac{1}{4}+(n+\frac{1}{2})\sum_{i}\frac{\varepsilon_{i}}{b_{i}}+\frac{1}{4}(\sum_{i}\frac{\varepsilon_{i}}{b_{i}})^{2})}(t^{\frac{\varepsilon+2n+1}{2}}+t^{\frac{-(\varepsilon+2n+1)}{2}}).\] 
Following \cite{gukovmanolescu}, fix $\varepsilon_{2}$ and $\varepsilon_{3}$ and split into two cases based on the value of $\varepsilon_{1}$. If $\varepsilon_{1}=-1$, observe that $b_{1}b_{2}b_{3}(1+\sum_{i}\frac{\varepsilon_{i}}{b_{i}})= \alpha_k$ for some $k \in \{1,2,3,4\}$. The corresponding summation over $n$ for this triple of $\varepsilon$'s is 
\begin{equation} \label{k sum} -\varepsilon_{2}\varepsilon_{3}\sum_{n\geq 0}q^{pn^{2}+\alpha_{k}n+\frac{\alpha_{k}^{2}}{4p}}(t^{\frac{\varepsilon_{2}+\varepsilon_{3}+2n}{2}}+t^{\frac{-(\varepsilon_{2}+\varepsilon_{3}+2n)}{2}}).\end{equation} 
On the other hand, when $\varepsilon_{1}=1$, we can replace $n$ with $n-1$ in the corresponding sum to get \begin{equation}\label{j sum} \varepsilon_{2}\varepsilon_{3}\sum_{n\geq 1}q^{pn^{2}-\alpha_{j}n+\frac{\alpha_{j}^{2}}{4p}}(t^{\frac{\varepsilon_{2}+\varepsilon_{3}+2n}{2}}+t^{\frac{-(\varepsilon_{2}+\varepsilon_{3}+2n)}{2}}),\end{equation} where for each $k$ the corresponding $j$ is given by

\begin{center}\begin{tabular}{c | c c c c}
$k$ & 1 & 2 & 3 & 4 \\ \hline
$j$ & 4 & 3 & 2 & 1.
\end{tabular}
\end{center}

\begin{remark}
In \cite{gukovmanolescu}, it is incorrectly claimed that $j=k$ for each $k$. This fact does not change the outcome of their calculations, but it does affect ours.
\end{remark} Summing over all four possible values of $(\varepsilon_{2},\varepsilon_{3})$ gives eight sums, each of which has exponent on $q$ of the form $pn^{2} \pm\alpha_{k}n+\frac{n^{2}}{4p}$ as in (\ref{k sum}) and (\ref{j sum}). The sums involving $+\alpha_{k}n$ begin at $n=0$ and the sums involving $-\alpha_{k}n$ begin at $n=1$. The four values of $(\varepsilon_{1},\varepsilon_{2},\varepsilon_{3})$ for which $\varepsilon_{2}=-\varepsilon_{3}$ contribute \[\sum_{n \geq 0}q^{pn^{2}+\alpha_{2}n+\frac{\alpha_{2}^{2}}{4p}}(t^n+t^{-n})-\sum_{n\geq 1}q^{pn^{2}-\alpha_{2}n+\frac{\alpha_{2}^{2}}{4p}}(t^n+t^{-n});\] \begin{equation}\label{alpha 3 sums} \sum_{n \geq 0}q^{pn^{2}+\alpha_{3}n+\frac{\alpha_{3}^{2}}{4p}}(t^n+t^{-n})-\sum_{n\geq 1}q^{pn^{2}-\alpha_{3}n+\frac{\alpha_{3}^{2}}{4p}}(t^n+t^{-n}),\end{equation} whereas when $\varepsilon_{2}=\varepsilon_{3}$ we have \[-\sum_{n \geq 0}q^{pn^{2}+\alpha_{4}n+\frac{\alpha_{4}^{2}}{4p}}(t^{n+1}+t^{-(n+1)})+\sum_{n\geq 1}q^{pn^{2}-\alpha_{4}n+\frac{\alpha_{4}^{2}}{4p}}(t^{n-1}+t^{-(n-1)});\] \begin{equation}\label{eqn:alpha1}-\sum_{n \geq 0}q^{pn^{2}+\alpha_{1}n+\frac{\alpha_{1}^{2}}{4p}}(t^{n-1}+t^{-(n-1)})+\sum_{n\geq 1}q^{pn^{2}-\alpha_{1}n+\frac{\alpha_{1}^{2}}{4p}}(t^{n+1}+t^{-(n+1)}).\end{equation} For $t=1$ and $\alpha_{k} \geq 0$, each of the above collapse to the false theta functions $\tilde{\Psi}^{(\alpha_{k})}_{p}$ into which $\Zhat$ is decomposed in \cite{gukovmanolescu}. The only case in which $\alpha_k < 0$ for some $k$ is $\Sigma(2,3,5),$ for which $\alpha_{1}=-1$. 
We momentarily postpone this case and take $(b_{1},b_{2},b_{3})\neq(2,3,5)$. Working with (\ref{alpha 3 sums}), we write $pn^2\pm n\alpha_3 +\frac{a_3^2}{4p} = p(n\pm\frac{\alpha_3}{2p})^2$ and perform the changes of variables $m = 2pn \pm \alpha_3$. This gives 
\[\sum_{\substack{m \geq 0 \\ m \equiv \alpha_3 \ (2p)}} q^{\frac{m^2}{4p}}(t^{\frac{m-\alpha_3}{2p}}+t^{-\frac{m-\alpha_3}{2p}}) -\sum_{\substack{m \geq 0 \\ m \equiv - \alpha_3 \ (2p)}} q^{\frac{m^2}{4p}}(t^{\frac{m+\alpha_3}{2p}}+t^{-\frac{m+\alpha_3}{2p}}).
\]
The calculation is the same when $\alpha_3$ is replaced with $\alpha_2$.
When $\varepsilon_{1}=\varepsilon_{3}=1$, we get the sums
\[-\sum_{\substack{m \geq 0 \\ m \equiv \alpha_4 \ (2p)}} q^{\frac{m^2}{4p}}(t^{\frac{m-\alpha_4}{2p}+1}+t^{-(\frac{m-\alpha_4}{2p}+1)}) + \sum_{\substack{m \geq 0 \\ m \equiv -\alpha_4 \ (2p)}} q^{\frac{m^2}{4p}}(t^{\frac{m+\alpha_4}{2p}-1}+t^{-(\frac{m+\alpha_4}{2p}-1)})
\] and when $\varepsilon_{2}=\varepsilon_{3}=-1$ we get
\[ -\sum_{\substack{m \geq 0 \\ m \equiv \alpha_1 \ (2p)}} q^{\frac{m^2}{4p}}(t^{\frac{m-\alpha_1}{2p}-1}+t^{-(\frac{m-\alpha_1}{2p}-1)}) + \sum_{\substack{m \geq 0 \\ m \equiv -\alpha_1 \ (2p)}} q^{\frac{m^2}{4p}}(t^{\frac{m+\alpha_1}{2p}+1}+t^{-(\frac{m+\alpha_1}{2p}+1)}).\]

If $(b_{1},b_{2},b_{3})\neq(2,3,5)$ we are done. We conclude with the special case of the Poincar\'e homology sphere. The argument is the same up through the calculation of (\ref{eqn:alpha1}).
In this case, we have that
\begin{align*} &- \sum_{n \geq 1} q^{30n^2-n+\frac{1}{120}}(t^{n-1}+t^{-(n-1)}) + \sum_{n \geq 0} q^{30n^2+n+\frac{1}{120}}(t^{n+1}+t^{-(n+1)}) \\ &= -\sum_{\substack{m \geq 0 \\ m \equiv -1 \ (60)}} q^{\frac{m^2}{120}}(t^{\frac{m-59}{60}}+t^{-(\frac{m-59}{60})}) + \sum_{\substack{m \geq 0 \\ m \equiv 1 \ (60)}} q^{\frac{m^2}{120}}(t^{\frac{m+59}{60}}+t^{-(\frac{m+59}{60})})
\end{align*}
and the bounds on the sums on the left hand side do not agree with those in (\ref{eqn:alpha1}). The solution is to subtract $2q^\frac{1}{120}(t+t^{-1})$ from (\ref{eqn:alpha1}), as they only disagree in the sign of their constant term.
\end{proof}

\subsection{Example Calculations} 
Below are calculations of $\Zhathat_{\Sigma}(\zeta, q)$ for various specializations of $\zeta$ for the Brieskorn spheres $\Sigma(2,3,5)$ and $\Sigma(2,7,15)$.
\small
\begin{center}
\begin{tabular}{|l|l|}
\hline
$\zeta$ & $\Zhathat_{\Sigma}(\zeta, q)$ \\
\hline
& \\
$1$ & $ 2q^{-3/2}- q^{-3/2}(1 + q + q^3 + q^7 - q^8 - q^{14} - q^{20} - q^{29} + q^{31} + q^{42} + q^{52} + \ldots)$ \\  & \\
$-1$ & $-2q^{-3/2}- q^{-3/2}(-1 + q + q^3 + q^7 + q^8 + q^{14} + q^{20} - q^{29} + q^{31} - q^{42}  + \ldots$)\\ & \\
$e^{\frac{2\pi i}{3}}$ & $ -q^{-3/2} - q^{-3/2}(-\frac{1}{2} + q + q^3 + q^7 + \frac{1}{2} q^8 + \frac{1}{2} q^{14} + \frac{1}{2} q^{20} - q^{29} - \frac{1}{2} q^{31}  + \ldots)$\\ & \\
$i$ & $-q^{-3/2}(q + q^3 + q^7 - q^{29} - q^{31} + q^{69} + q^{85} + q^{99} - q^{143} - q^{161} - q^{185}  +  \ldots)$\\ 
& \\
\hline
\end{tabular}
\end{center}
\begin{center}
\textsc{Table 1}: Examples of $\Zhathat_{\Sigma}(\zeta, q)$ for $\Sigma(2,3,5)$.
\end{center}

\small
\begin{center}
\begin{tabular}{|l|l|}
\hline
$\zeta$ & $\Zhathat_{\Sigma}(\zeta, q)$ \\
\hline
& \\
$1$ & $-q^{5/2}( -q^4 + q^9 + q^{17} - q^{26} + q^{87} - q^{106} - q^{130} + q^{153} - q^{275}+ q^{308}  + \dots) $\\  & \\
$-1$ & $-q^{5/2}(q^4 + q^9 + q^{17} + q^{26} + q^{87} + q^{106} + q^{130} + q^{153} - q^{275} - q^{308} + \dots) $ \\ & \\
$e^{\frac{2\pi i}{3}}$ & $-q^{5/2}(\frac{1}{2}q^4 + q^9 + q^{17} + \frac{1}{2}q^{26} + \frac{1}{2}q^{30} - \frac{1}{2}q^{153} - q^{275} - \frac{1}{2}q^{308} - \frac{1}{2}q^{348} + \dots) $
\\ & \\
$i$ & $-q^{5/2}( q^9 + q^{17} + q^{87} - q^{153} - q^{275} + q^{385} + q^{615} + q^{671} - q^{1027} - q^{1099} + \dots) $
\\ 
& \\
\hline
\end{tabular}
\end{center}
\begin{center}
\textsc{Table 2}: Examples of $\Zhathat_{\Sigma}(\zeta, q)$ for $\Sigma(2,7,15)$.
\end{center}
\normalsize

In the above examples, we factor out a rational power of $q$ so that all other powers are integral. This can be done in general, and is explicitly realized for Brieskorn spheres in the following lemma:

\begin{lemma}\label{lem: congruence class}
Let $\alpha_k$, $1 \leq k \leq 4$, be as above. Then $\alpha_1^2 \equiv \alpha_2^2 \equiv \alpha_3^2 \equiv \alpha_4^2 \bmod{4p}$.
\end{lemma}
We will denote this common congruence class mod $4p$ by $w$ for the remainder of the paper.

\section{Radial Limits at Roots of Unity}\label{sec: radial limits}

In this section, we analyze the arithmetic properties of the coefficients $\varphi(n;t)$ which will ultimately allow for the calculation of radial limits at roots of unity in terms of particular $L$-functions. We first check that the coefficients of $\varphi(n;t)$ have the necessary properties for our method of calculation.

\begin{lemma} If $\zeta$ is a $j$th root of unity, then $\varphi(n; \zeta)$ is $2pj$-periodic and has mean value zero.
\end{lemma}
In order to calculate these radial limits, we will make use of the following general proposition. 

\begin{prop} \label{prop: mellin limits}
Let $C \colon \mathbb{Z} \to \mathbb{C}$ be a periodic function with mean value zero. Then the associated $L$-series $L(s, C) := \sum_{n \geq 1}\frac{C(n)}{n^s}, \ \Re(s) > 1$, extends holomorphically to all of $\mathbb{C}$ and the function $\sum_{n \geq 1} C(n)e^{-n^2t}, \ t > 0$, has the asymptotic expansion
\[\sum_{n \geq 1} C(n)e^{-n^2t} \sim \sum_{r \geq 0} L(-2r, C) \frac{(-t)^r}{r!}
\]
as $t \searrow 0$. Then numbers $L(-r,C)$ are given explicitly by 
\[L(-r, C) = -\frac{M^r}{r+1} \sum_{n=1}^M C(n) B_{r+1}\left(\frac{n}{M}\right) \quad (r=0, 1, \dots)
\]
where $B_k(x)$ is the $k$th Bernoulli polynomial and $M$ is any period of the function $C(n)$. 
\end{prop}

For details, see e.g. \cite{lawrencezagier} p. 98.

\subsection{Proof of Theorem \ref{thm: radial limits}}
Let $\xi$ be a root of unity and set $C(n):=\varphi(n;\zeta)\xi^\frac{n^{2}}{4p}$.
If $K$ is a period of $\xi$, then $C(n)$ is $2pjK$-periodic and has mean value zero since
$
C(2pjK-n)
=-C(n).$
Let \begin{equation}A_\zeta(q) := \sum_{n \geq 0} \varphi(n; \zeta)q^\frac{n^2}{4p}
\end{equation} 
and observe that 
\begin{align*}
A_\zeta(\xi e^{-t})=\sum_{n=1}^{\infty}C(n)e^{-n^2(t/4p)}.
\end{align*} By the previous proposition, the above has asymptotic expansion 
\[\sum_{r=0}^{\infty}L(-2r,C)\frac{(-t/4p)^r}{r!}\] as $t\searrow 0$ and limiting value \[A_\zeta(\xi):=\lim_{t\searrow 0}A_\zeta(\xi e^{-t})= L(0,C).\]
The analytic continuation of this $L$-function to $s=0$ is given by the sum
\[ -\sum_{n=1}^{2pjK} \left(\frac{n}{2pjK}-\frac{1}{2} \right) \varphi(n;\zeta)\xi^\frac{n^2}{4p}.
\]
Accounting for $C$ and $q^\Delta$ evaluated at $(\zeta, \xi)$ gives the desired formula. $\Box$ 

\section{Modular and Quantum Modular Forms}\label{sec: modularity discussion}

We begin with a brief introduction to the theory of modular forms of half-integral weight. For a more thorough discussion, see \cite{web, shimura}.
Let $\gamma = \begin{psmallmatrix} a & b \\ c & d \end{psmallmatrix} \in \text{SL}_2(\mathbb{Z})$ act on $\mathbb{H}$ by the linear fractional transformation
\[\gamma\tau := \frac{a \tau + b}{c \tau +d}.\]
We will be interested in the action of particular congruence subgroups $\Gamma$ of $\text{SL}_2(\mathbb{Z})$. Define 
\begin{align*}
\Gamma_1(N)&:= \left\{\begin{pmatrix} a & b \\ c & d
\end{pmatrix} \in \text{SL}_2(\mathbb{Z}) : \ a \equiv d \equiv 1 \bmod{N}, \ c \equiv 0 \bmod{N}\right\}; \\ \Gamma(N)&:= \left\{\begin{pmatrix} a & b \\ c & d
\end{pmatrix} \in \text{SL}_2(\mathbb{Z}) : \ a \equiv d \equiv 1 \bmod{N}, \ b \equiv c \equiv 0 \bmod{N}\right\}.
\end{align*}
The above are congruence subgroups of $\text{SL}_2(\mathbb{Z})$ of level $N$. Note that $\Gamma(N) \subseteq \Gamma_1(N)$. The equivalence classes in $\mathbb{P}^1(\mathbb{Q}) = \mathbb{Q} \cup \{ i \infty \}$ under the action of a congruence subgroup $\Gamma$ are called the \textit{cusps} of $\Gamma$. 

To state the appropriate transformation law for half-integral weight modular forms, we need the following definitions. For odd $d$, define 
\[\varepsilon_d := \begin{cases}
1 & \text{if } d \equiv 1 \bmod{4}; \\
i & \text{if } d \equiv 3 \bmod{4},
\end{cases}
\]
and let $\left(\frac{\cdot}{\cdot}\right)$ denote the Jacobi symbol. Throughout, we let $\sqrt{z}$ be the branch of the square root with argument in $(-\pi/2, \pi/2]$.
For functions $f\colon \mathbb{H} \to \mathbb{C}$, the Petersson slash operator of weight $k \in \frac{1}{2}\mathbb{Z}$ for $\gamma = \begin{psmallmatrix} a & b \\ c & d \end{psmallmatrix} \in \text{SL}_2(\mathbb{Z})$ is defined by
\[
f|_k\gamma(\tau):= \begin{cases}
(c\tau + d)^{-k} f(\gamma \tau)& \text{ if } k \in \mathbb{Z}; \\
\varepsilon_d^{2k}\left(\frac{c}{d}\right)(c\tau + d)^{-k}f(\gamma \tau) & \text{ if } k \in \frac{1}{2} + \mathbb{Z},
\end{cases}\]
where additionally one must require that $\gamma$ is contained in a congruence subgroup of level 4 when $k$ is not an integer. We can now state the following:

\begin{definition}
Let $\Gamma$ be a congruence subgroup of level $N$ such that $4 \mid N$. We say that a holomorphic function $f \colon \mathbb{H} \to \mathbb{C}$ is a modular form (resp. cusp form) of weight $k \in \frac{1}{2}\mathbb{Z}$ with multiplier $\chi$ for $\Gamma$ if:
\begin{enumerate}
\item for all $\gamma \in \Gamma$, the function $f$ satisfies $f -\overline{\chi}(\gamma) f|_{k}\gamma=0$, and
\item for all $\gamma \in \emph{SL}_2(\mathbb{Z})$, $(c\tau + d)^{-k}f( \gamma \tau)$ is bounded (resp. vanishes) as $\tau \to i\infty$. 
\end{enumerate}
\end{definition}

While our work will make contact with modular forms as described above, the modular objects of primary interest will be quantum modular forms. This term, coined by Zagier in 2010, was inspired in part by the examples arising from quantum field theory and quantum invariants of $3$-manifolds such as the WRT invariants \cite{zagier}.

Fix a congruence subgroup $\Gamma$ of level $N$ such that $4 \mid N$, and suppose $\mathcal{Q} = \mathbb{Q}\backslash S$ where $S$ is discrete and $\mathcal{Q}$ is closed under the action of $\Gamma$.
We define a \textit{quantum modular form} of weight $k$ with multiplier $\chi$ for $\Gamma$ to be a function $f\colon \mathcal{Q} \to \mathbb{C}$ such that for all $\gamma = \begin{psmallmatrix} a & b \\ c & d
\end{psmallmatrix} \in \Gamma$, the functions $h_\gamma\colon \mathcal{Q}\backslash \{\gamma^{-1}(i \infty)\}) \to \mathbb{C}$,
\begin{equation}h_\gamma(x) := f(x) - \overline{\chi}(\gamma)f|_k \gamma(x)\end{equation}
extends to some ``nice" function on $\mathbb{R}$. The set $\mathcal{Q}$ is called the quantum set of $f$.

\begin{remark} This definition is intentionally vague, as it was built to fit the particular examples naturally arising from disparate areas of study. While this definition is still under construction, there have been a variety of alterations made to the original definition that are still considered to fall under the ``quantum modular" umbrella.
\end{remark}

Ultimately, we will need the following lemma in order to renormalize the powers of $q$ that we encounter in Section \ref{sec: modularity proof}. 
Using the definitions, one can verify the following:
\begin{lemma}\label{lem:mult by j}
If $\psi(\tau)$ is a quantum modular form of weight $1/2$ for $\Gamma(4pj)$ with multiplier $\chi$, then $\psi(j\tau)$ is a quantum modular form of weight $1/2$ for $\Gamma(4pj^2)$ with the same multiplier.
\end{lemma}

\subsection{Eichler Integrals}\label{sec: eichler}

The utilization of Eichler integrals to construct quantum modular forms has its roots in the work of Lawrence and Zagier previously discussed \cite{lawrencezagier}. Many authors have since extended and generalized this procedure to systematically construct families of quantum modular forms: see e.g. \cite{folsomonorhoades, bringmannrolen, goswamiosburn}. 
Here we sketch the procedure by which the authors of \cite{bringmannrolen} construct quantum modular forms, modified to fit our context. For references that reflect these arguments with our particular congruence subgroups and multiplier systems, see e.g. \cite{bringmannmilaswalgebras, goswamiosburn}. 

Suppose a function $F(\tau)$ for $\tau \in \mathbb{H}$ may be written as
\[F(\tau) = \sum_{n \geq 0} a(n) q^\frac{n^2}{4pj}, \quad \quad (q=e^{2\pi i \tau})
\]
for some integers $p, \ j$. Further suppose that
\[f(\tau):= \sum_{n \geq 0} n a(n) q^\frac{n^2}{4pj}
\]
is a cusp form of weight $3/2$ for $\Gamma_1(N)$.  
We consider the \textit{non-holomorphic Eichler integral} of $f$, given by 
\[F^*(\tau):= \int_{\overline{\tau}}^{i\infty} \frac{f(\omega)}{\sqrt{-i(\omega-\tau)}} \ d\omega, \quad \quad (\tau \in \mathbb{H}^-).
\]
Bringmann and Rolen show that, after suitable renormalization, the functions $F(\tau)$ and $F^*(\tau)$ ``agree to infinite order" at any $x \in \mathbb{Q}$. That is, for any $x$ there exists a sequence $\beta(0), \beta(1), \dots$ such that as $t \to 0^+$, 
\[F\left(x+\frac{it}{2\pi}\right) \sim \sum_{r \geq 0} \beta(r) \frac{(-t)^r}{r!} \ \text{ and } \ F^*\left(x-\frac{it}{2\pi}\right) \sim \sum_{r \geq 0} \beta(r) \frac{t^r}{r!}.\]
This is accomplished by first integrating $F^*$ term-by-term to obtain a series expansion for $F^*(\tau)$ involving $\Gamma$-factors. Then using Proposition \ref{prop: mellin limits} and more general tools for studying the Mellin transform of error functions, they obtain the asymptotic expansions of both of these functions and verify that they agree in the above sense.

The function $F^*(\tau)$ admits an explicit obstruction to modularity from its definition; for $\tau \in \mathbb{H}^-$ and $\gamma = \begin{psmallmatrix} a & b \\ c & d \end{psmallmatrix} \in \Gamma_1(N)$, we have
\[F^*(\tau) - \textstyle{\left(\frac{-4}{d}\right)}F^*|_k\gamma \left(\tau\right)= r_\gamma(\tau),
\]
where 
\[r_{\gamma}(\tau):= \int_{\gamma^{-1}(i\infty)}^{i \infty}\frac{f(\omega)}{\sqrt{-i(\omega-\tau)}} \ d\omega.
\]
which extends to a $C^\infty$ function on $\partial \mathbb{H}^- = \mathbb{R}$ which is real-analytic on $\mathbb{R}\backslash \{\gamma^{-1}(i\infty)\}$
and gives $h_\gamma$ for the resulting quantum modular form.

\section{Proofs of Main Theorems \label{sec: modularity proof}}
Below we prove two results regarding the quantum modularity of $\Zhathat(t,q)$.
\subsection{Proof of Theorem \ref{thm:modularity}}
In light of the work summarized in the preceding section, it suffices to show that \begin{equation} \label{eqn: series2analyze} \sum_{n \geq 0} n \varphi(n;\zeta) q^\frac{n^2}{4pj}\end{equation}
is a cusp form for $\Gamma(4pj)$.
Then the results of Section \ref{sec: eichler} used in conjunction with Lemma \ref{lem:mult by j} will imply 
\[\sum_{n \geq 0} \varphi(n;\zeta) q^\frac{n^2}{4p}\]
is a quantum modular form.
We begin with an elementary lemma which will be useful for simplifying our expressions later. 
\begin{lemma} \label{lem: congruence mod 4pj}
Let $0 \leq n < 2pj$ be such that $n \equiv \pm \alpha_k \bmod{2p}$ for some $k$. Then we have $n^2 \equiv w + 4pi \bmod{4pj}$ for some $0 \leq i < j$, where $w$ is the common congruence class mod $4p$ of the $\alpha_k^2$'s coming from Lemma $\ref{lem: congruence class}$. 
\end{lemma}
We are now ready to analyze (\ref{eqn: series2analyze}).
Since $\varphi(n;\zeta)$ is $2pj$-periodic, we have
\begin{align*}
\sum_{n \geq 0} n \varphi(n;\zeta) q^\frac{n^2}{4pj} 
&= \sum_{0 \leq \alpha < 2pj} \varphi(\alpha; \zeta) \sum_{n \geq 0} (2pjn+\alpha)q^\frac{(2pjn+\alpha)^2}{4pj}.
\end{align*}
Every $\alpha$ for which $\varphi(\alpha;\zeta)$ is nonzero satisfies $\alpha \equiv \pm \alpha_k \bmod{2p}$ for some $k$.
Thus, we can write this sum as
\begin{align*}  & \sum_{\substack{0 \leq \alpha < 2pj \\ \alpha \equiv \alpha_k (2p)}} \varphi(\alpha, \zeta)  \sum_{n \geq 0} (2pjn+\alpha)q^\frac{(2pjn+\alpha)^2}{4pj} \\ + &\sum_{\substack{0 < \alpha \leq 2pj \\ \alpha \equiv \alpha_k (2p)}} \varphi(2pj - \alpha, \zeta) \sum_{n \geq 0} (2pjn+(2pj-\alpha))q^\frac{(2pjn+(2pj-\alpha))^2}{4pj}.\end{align*}
Using the fact that $\varphi(n;\zeta)$ is odd and $2pj$-periodic, the second set of sums can be rewritten as
\[\sum_{\substack{0 < \alpha \leq 2pj \\ \alpha \equiv \alpha_k (2p)}} -\varphi(\alpha, \zeta) \sum_{n \geq 0} (2pj(n+1)-\alpha)q^\frac{(2pj(n+1)-\alpha))^2}{4pj}.\]
Reindexing by $n+1 \mapsto -n$ and combining with the first set of sums gives
\begin{equation}\label{eqn: main sum} \sum_{n \geq 0} n \varphi(n;\zeta) q^\frac{n^2}{4pj} = \sum_{\substack{0 \leq \alpha < 2pj \\ \alpha \equiv \alpha_k (2p)}} \varphi(\alpha; \zeta) \sum_{\substack{n \in \mathbb{Z}\\ n \equiv \alpha (2pj)}} nq^\frac{n^2}{4pj}.
\end{equation}
The inner sum of the above equation is a theta function which is modular of weight $3/2$. More precisely, define \[\Theta(\tau; k, M):= \sum_{\substack{n \in \mathbb{Z}\\ n \equiv k (M)}} nq^\frac{n^2}{2M}.
\]
By Proposition 2.1 of \cite{shimura}, we have that for $\gamma \in \Gamma_1(2M)$ that
\[\Theta(\gamma \tau; k, M) = e^\frac{\pi i a b k^2}{M}\varepsilon_d^{-3} \left( \frac{2Mc}{d}\right) (c\tau + d)^{3/2} \Theta(\tau; ak, M),\]
and since $k \equiv ak \bmod{M}$, we have
\[\Theta(z; ak, M) = \Theta(z; k, M).
\]
By Lemma 6.1, every $n$ for which the coefficient of $q^\frac{n^2}{4pj}$ in (\ref{eqn: main sum}) is nonzero satisfies $n^2 \equiv w + 4pi \bmod{4pj}$ for some $0 \leq i < j$. 
Then
\[e^\frac{\pi i a b n^2}{2pj} = e^\frac{\pi i a b (w + 4pi)}{2pj}
\]
for some $0 \leq i < j$. Then one may group the $n$'s based on the corresponding $i$ to get $j$ functions $f_{i}(\tau)$ which for $\gamma = \begin{psmallmatrix} a & b \\ c & d \end{psmallmatrix} \in \Gamma_1(4pj)$ satisfy
\[f_{i}(\gamma \tau) = e^\frac{\pi i a b (w + 4pi)}{2pj} \varepsilon_d^{-3} \left( \frac{4pjc}{d}\right)(c\tau + d)^{3/2} f_{i}(\tau).
\]
Note the dependence of this transformation law on $w$. If one restricts to $\gamma \in \Gamma(4pj) \subset \Gamma_1(4pj)$, the multipliers for each $f_i$ become identical. Thus the sum of the $f_i$'s transform together as a cusp form on $\Gamma(4pj)$.
$\Box$

\subsection{Proof of Theorem \ref{thm:diffmodular}}

As in the study of the Jacobi Triple Product formula, one is often able to generate a modular object of dual-weight by differentiating with respect to one variable (see e.g. \cite{coogan-ono}). Following this approach, we find a second infinite family of quantum invariants by differentiating $\Zhathat_\Sigma(t,q)$, summand by summand, with respect to $t$. Our contribution to this principle is Theorem \ref{thm:diffmodular}. Here we offer of proof of this result.

Fix $\zeta$ a $j$th root of unity. Consider the series
\[A'_\zeta(\tau):= \sum_{n \geq 0} \varphi'(n;\zeta) q^\frac{n^2}{4p},
\]
where $\varphi'(n; \zeta)$ is the derivative of $\varphi(n; t)$ evaluated at $t = \zeta$. By Theorem \ref{thm: q-series}, this is
\[\varphi'(n;\zeta) := \begin{cases}
\frac{n \mp (\alpha_1+2p)}{4p}(\zeta^{\frac{\mp n + (\alpha_1+2p)}{2p}} - \zeta^{\frac{\pm n - (\alpha_1+2p)}{2p}}) & n \equiv \pm \alpha_1 \bmod{2p},  \\
- \frac{n\mp \alpha_k}{4p}(\zeta^\frac{\mp n + \alpha_k}{2p} - \zeta^{\frac{\pm n - \alpha_k}{2p}}) & n \equiv \pm \alpha_k \bmod{2p}, \ k=2,3\\
 \frac{n \mp (\alpha_4 -2p)}{4p}(\zeta^{\frac{\mp n + (\alpha_4-2p)}{2p}} - \zeta^{\frac{\pm n - (\alpha_4-2p)}{2p}}) & n \equiv \pm \alpha_4 \bmod{2p}, \\ 
0 & otherwise. 
\end{cases}\]
Note that we may write $A'_\zeta(\tau)$ as
\[ \sum_{n \geq 0} n \psi(n;\zeta)q^\frac{n^2}{4p} +  \sum_{n \geq 0} \chi(n;\zeta)q^\frac{n^2}{4p}, 
\]
where
\begin{align*}\psi(n;\zeta) &:= \begin{cases}
\frac{1}{4p}(\zeta^{\frac{\mp n + (\alpha_1+2p)}{2p}} - \zeta^{\frac{\pm n - (\alpha_1+2p)}{2p}}) & n \equiv \pm \alpha_1 \bmod{2p},  \\
- \frac{1}{4p}(\zeta^\frac{\mp n + \alpha_k}{2p} - \zeta^{\frac{\pm n - \alpha_k}{2p}}) & n \equiv \pm \alpha_k \bmod{2p}, \ k=2,3\\
\frac{1}{4p}(\zeta^{\frac{\mp n + (\alpha_4-2p)}{2p}} - \zeta^{\frac{\pm n - (\alpha_4-2p)}{2p}}) & n \equiv \pm \alpha_4 \bmod{2p}, \\ 
0 & otherwise;
\end{cases} \\
\chi(n;\zeta) &:= \begin{cases}
\mp \frac{(\alpha_1+2p)}{4p}(\zeta^{\frac{\mp n + (\alpha_1+2p)}{2p}} - \zeta^{\frac{\pm n - (\alpha_1+2p)}{2p}}) & n \equiv \pm \alpha_1 \bmod{2p},  \\
\pm \frac{\alpha_k}{4p}(\zeta^\frac{\mp n + \alpha_k}{2p} - \zeta^{\frac{\pm n - \alpha_k}{2p}}) & n \equiv \pm \alpha_k \bmod{2p}, \ k=2,3\\
 \mp \frac{(\alpha_4 -2p)}{4p}(\zeta^{\frac{\mp n + (\alpha_4-2p)}{2p}} - \zeta^{\frac{\pm n - (\alpha_4-2p)}{2p}}) & n \equiv \pm \alpha_4 \bmod{2p}, \\ 
0 & otherwise. 
\end{cases}\end{align*}
Then $\psi(n;\zeta)$ is even and $2pj$-periodic and $\chi(n;\zeta)$ is odd and $2pj$-periodic. 
Following the same style of argument as Theorem \ref{thm:modularity}, one concludes that $\displaystyle{\sum_{n \geq 0} \chi(n;\zeta)q^\frac{n^2}{4p}}$ is a quantum modular form of weight $1/2$ on $\Gamma(4pj^2)$. 

To analyze $\displaystyle{\sum_{n \geq 0} n \psi(n;\zeta)q^\frac{n^2}{4p}}$, note that the series $\displaystyle{\sum_{n \geq 0} \psi(n;\zeta)q^\frac{n^2}{4p}}$ is modular but may not be a cusp form. This requires us to appeal to a more general result of Goswami and Osburn (Theorem 1.1 of \cite{goswamiosburn}) which gives a careful treatment of this more general case. Their result tells us that 
\[\sum_{n \geq 0} n \psi(n;\zeta)q^\frac{n^2}{4pj}
\]
is a quantum modular form on $Q_{2pj}$ with respect to $\Gamma_1(4pj)$, where
\[Q_{2pj} := \{x \in \mathbb{Q} \ : \ x \text{ is } \Gamma_1(4pj)\text{-equivalent to } i\infty\}.
\]
Note that one must still utilize Lemmas \ref{lem:mult by j} and \ref{lem: congruence mod 4pj} in order to contend with the supports of these series.
This ultimately allows us to conclude that $\sum_{n \geq 0} n \psi(n;\zeta)q^\frac{n^2}{4p}$ is a quantum modular form of weight $3/2$ as desired.
$\Box$

\bibliographystyle{amsplain}
\bibliography{Biblio.bib}

\end{document}